\documentclass[11pt]{article}
\usepackage{latexsym,amssymb,amsmath,enumerate,geometry,float,cite,tikz,setspace}
\newtheorem{definition}{Definition}
\newtheorem{theorem}[definition]{Theorem}
\newtheorem{corollary}[definition]{Corollary}

\usepackage{lineno}

\begin{document}


\onehalfspace

\title{Relating $2$-Rainbow Domination to Roman domination}

\author{Jos\'{e} D. Alvarado$^1$, Simone Dantas$^1$, and Dieter Rautenbach$^2$}

\date{}

\maketitle

\begin{center}
{\small 
$^1$ Instituto de Matem\'{a}tica e Estat\'{i}stica, Universidade Federal Fluminense, Niter\'{o}i, Brazil\\
\texttt{josealvarado.mat17@gmail.com, sdantas@im.uff.br}\\[3mm]
$^2$ Institute of Optimization and Operations Research, Ulm University, Ulm, Germany\\
\texttt{dieter.rautenbach@uni-ulm.de}
}
\end{center}

\begin{abstract}
For a graph $G$, 
let $\gamma_R(G)$ and $\gamma_{r2}(G)$ denote 
the Roman domination number of $G$ and 
the $2$-rainbow domination number of $G$,
respectively.
It is known that $\gamma_{r2}(G)\leq \gamma_R(G)\leq \frac{3}{2}\gamma_{r2}(G)$.
Fujita and Furuya (Difference between 2-rainbow domination and Roman domination in graphs, Discrete Applied Mathematics 161 (2013) 806-812) present some kind of characterization 
of the graphs $G$ for which $\gamma_R(G)-\gamma_{r2}(G)=k$
for some integer $k$.
Unfortunately, their result does not lead to an algorithm that allows to recognize these graphs efficiently.

We show that for every fixed non-negative integer $k$,
the recognition 
of the connected $K_4$-free graphs $G$ with $\gamma_R(G)-\gamma_{r2}(G)=k$
is NP-hard,
which implies that there is most likely no good characterization of these graphs.
We characterize the graphs $G$ such that $\gamma_{r2}(H)=\gamma_R(H)$ for every induced subgraph $H$ of $G$,
and collect several properties of the graphs $G$ with $\gamma_R(G)=\frac{3}{2}\gamma_{r2}(G)$.
\end{abstract}

{\small 

\medskip


\medskip


}

\pagebreak

\section{Introduction}\label{section1}

We consider finite, simple, and undirected graphs, and use standard terminology and notation. 

Rainbow domination was introduced in \cite{bhr}.
Here we consider the special case of $2$-rainbow domination. 
A {\it $2$-rainbow dominating function} of a graph $G$ is a function $f:V(G)\to 2^{\{ 1,2\}}$ 
such that $\bigcup_{v\in N_G(u)}f(v)=\{ 1,2\}$ for every vertex $u$ of $G$ with $f(u)=\emptyset$.
For a set $X$ of vertices of $G$, let $|f(X)|=\sum_{u\in X}|f(u)|$, and
let the {\it weight $w(f)$} of $f$ be $|f(V(G))|$.
The {\it $2$-rainbow domination number $\gamma_{r2}(G)$} of $G$ 
is the minimum weight of a $2$-rainbow dominating function of $G$,
and a $2$-rainbow dominating function of $G$ of weight $\gamma_{r2}(G)$ is {\it minimum}.

Roman domination was introduced in \cite{s}.
A {\it Roman dominating function} of a graph $G$ is a function $g:V(G)\to \{ 0,1,2\}$ 
such that every vertex $u$ of $G$ with $g(u)=0$ has a neighbor $v$ with $g(v)=2$.
For a set $X$ of vertices of $G$, let $g(X)=\sum_{u\in X}g(u)$, and
let the {\it weight $w(g)$} of $g$ be $g(V(G))$.
The {\it Roman domination number $\gamma_R(G)$} of $G$
is the minimum weight of a Roman dominating function of $G$,
and a Roman dominating function of $G$ of weight $\gamma_R(G)$ is {\it minimum}.

For a positive integer $k$, let $[k]$ be the set of positive integers at most $k$.

\medskip

\noindent The definitions of $2$-rainbow domination on the one hand and Roman domination on the other hand clearly exhibit certain similarities.
It is therefore not surprising that these notions are related.
For later reference, we include the simple proof of the following known results.

\begin{theorem}[Wu and Xing \cite{wx}, Chellali and Rad \cite{cr}, and Fujita and Furuya \cite{ff}]\label{theorembefore}
If $G$ is a graph, then $\gamma_{r2}(G)\leq \gamma_R(G)\leq \frac{3}{2}\gamma_{r2}(G)$.
\end{theorem}
{\it Proof:} These inequalities follow immediately from two simple observations:
If $g$ is a Roman dominating function of $G$, then 
$$f:V(G)\to 2^{\{ 1,2\}}: x\mapsto
\begin{cases}
\emptyset & \mbox{, if $g(x)=0$,}\\
\{ 1\} & \mbox{, if $g(x)=1$, and}\\
\{ 1,2\} & \mbox{, if $g(x)=2$}
\end{cases}
$$
is a $2$-rainbow dominating function of $G$ of weight $w(f)\leq w(g)$.
Similarly, if $f$ is a $2$-rainbow dominating function of $G$, 
and $|f^{-1}(\{ 1\})|\geq |f^{-1}(\{ 2\})|$, then
$$g:V(G)\to \{ 0,1,2\}: x\mapsto
\begin{cases}
0 & \mbox{, if $f(x)=\emptyset$,}\\
1 & \mbox{, if $f(x)=\{ 1\}$, and}\\
2 & \mbox{, otherwise}
\end{cases}
$$
is a Roman dominating function of $G$ of weight
\begin{eqnarray}
w(g) & = & |f^{-1}(\{ 1\})|+2|f^{-1}(\{ 2\})|+2|f^{-1}(\{ 1,2\})|\nonumber\\
& \leq & \frac{3}{2}|f^{-1}(\{ 1\})|+\frac{3}{2}|f^{-1}(\{ 2\})|+2|f^{-1}(\{ 1,2\})|\label{e1}\\
& \leq & \frac{3}{2}|f^{-1}(\{ 1\})|+\frac{3}{2}|f^{-1}(\{ 2\})|+3|f^{-1}(\{ 1,2\})|\label{e2}\\
&=& \frac{3}{2}w(f)\nonumber.
\end{eqnarray}
$\Box$

\medskip

\noindent Fujita and Furuya \cite{ff} present some kind of characterization 
of the connected graphs $G$ for which $\gamma_R(G)-\gamma_{r2}(G)=k$
for some non-negative integer $k$ at most $\frac{1}{2}\gamma_{r2}(G)$ (cf. Corollary 3.6 in \cite{ff}). 
Unfortunately, their result does not lead to an algorithm that allows to recognize these graphs efficiently.

In the present note we show that for every fixed non-negative integer $k$,
the recognition 
of the connected $K_4$-free graphs $G$ with $\gamma_R(G)-\gamma_{r2}(G)=k$
is NP-hard,
which implies that there is most likely no good characterization of these graphs.
In view of this negative result, 
we characterize the graphs $G$ such that $\gamma_{r2}(H)=\gamma_R(H)$ for every induced subgraph $H$ of $G$, 
and also establish a similar result for the equality $\gamma_R(H)=\frac{3}{2}\gamma_{r2}(H)$.
The graphs $G$ that satisfy $\gamma_R(G)=\frac{3}{2}\gamma_{r2}(G)$ seem far more restricted
and we collect several of their properties.

For further related results on these parameters refer to \cite{wr,adr1,adr2}.

\section{Results}

We begin with our hardness results.

\begin{theorem}\label{theorem1}
It is NP-hard to decide whether $\gamma_{r2}(G)=\gamma_R(G)$ for a given connected $K_4$-free graph $G$. 
\end{theorem}
{\it Proof:} We describe a reduction from {\sc 3Sat}. 
Therefore, let ${\cal F}$ be an instance of {\sc 3Sat} with $m$ clauses $C_1,\ldots,C_m$ over $n$ boolean variables $x_1,\ldots,x_n$. 
Clearly, we may assume that $m\geq 2$. 
We will construct a connected $K_4$-free graph $G$ whose order is polynomially bounded in terms of $n$ and $m$ 
such that ${\cal F}$ is satisfiable if and only if $\gamma_{r2}(G)=\gamma_R(G)$.  
\begin{itemize}
\item For every variable $x_i$, we create a copy $G_i$ of the diamond $K_4-e$, 
and denote the two vertices of degree $3$ in $G_i$ by $x_i$ and $\bar{x}_i$. 
\item For every clause $C_j$, we create a vertex $C_j$.
\item For every literal $x\in\{ x_1,\ldots,x_n \}\cup\{ \bar{x}_1,\ldots,\bar{x}_n\}$ and every clause $C_j$ such that $x$ appears in $C_j$, 
we add the edge $xC_j$. 
\item Finally, we add an induced path $uvw$ of order $3$ and all possible edges between $\{ u,w\}$ and $\{ C_1,\ldots, C_m\}$. 
\end{itemize}
This completes the construction of $G$. 
Clearly, $G$ is connected and $K_4$-free, and has order $4n+m+3$. 

Let $f$ be a $2$-rainbow dominating function of $G$.
It is easy to see that $|f(V(G_i))|\geq 2$ for every $i\in [n]$, and $|f(\{ C_1,\ldots,C_m\}\cup \{ u,v,w\})|\geq 2$,
which implies $\gamma_{r2}(G)\geq 2n+2$.
Since 
$$x\mapsto
\begin{cases}
\{ 1\} &\mbox{, $x\in \{ u,x_1,\ldots,x_n\}$,}\\
\{ 2\} & \mbox{, $x\in \{ w,\bar{x}_1,\ldots,\bar{x}_n\}$, and}\\
\emptyset & \mbox{, otherwise}
\end{cases}
$$
defines a $2$-rainbow dominating function of weight $2n+2$,
we obtain $\gamma_{r2}(G)=2n+2$,
which implies $\gamma_R(G)\geq 2n+2$.

In remains to show that ${\cal F}$ is satisfiable if and only if $\gamma_R(G)=2n+2$.

Suppose that $\gamma_R(G)=2n+2$.
Let $g$ be a minimum Roman dominating function of $G$.
It is easy to see that $g(V(G_i))\geq 2$ for every $i\in [n]$, and $g(\{ C_1,\ldots,C_m\}\cup \{ u,v,w\})\geq 2$.
Since $\gamma_R(G)=2n+2$, all these inequalities are satisfied with equality,
and considering $u$, $v$, and $w$, it follows easily that $g(v)=2$, and that every vertex $C_j$
has a neighbor $x$ in $\{ x_1,\ldots,x_n\}\cup \{ \bar{x}_1,\ldots,\bar{x}_n\}$ with $g(x)=2$.
Therefore, these latter vertices indicate a satisfying truth assignment for ${\cal F}$.

Conversely, suppose that ${\cal F}$ is satisfiable, and consider a satisfying truth assignment.
The function
$$x\mapsto
\begin{cases}
2 &\mbox{, $x=v$},\\
2 &\mbox{, $x\in \{ x_1,\ldots,x_n\}\cup \{ \bar{x}_1,\ldots,\bar{x}_n\}$ such that $x$ is true, and}\\
0 & \mbox{, otherwise}
\end{cases}
$$
defines a Roman dominating function of $G$ of weight $2n+2$,
which implies $\gamma_R(G)=2n+2$, and completes the proof. $\Box$

\medskip

\noindent If $G$ is a graph, then 
$\gamma_{r2}(G\cup C_4)=\gamma_{r2}(G)+2$
and
$\gamma_R(G\cup C_4)=\gamma_R(G)+3$.
Furthermore, if $G$ has $k$ components $G_1,\ldots,G_k$,
the star $K_{1,k+2}$ has endvertices $u_1,\ldots,u_{k+2}$, and 
$G'$ arises from $G\cup K_{1,k+2}$ by adding an edge between $u_i$ and one vertex of $G_i$ for every $i\in [k]$,
then $G'$ is connected, and satisfies
$\gamma_{r2}(G')=\gamma_{r2}(G)+2$
and
$\gamma_R(G')=\gamma_R(G)+2$.
In combination with Theorem \ref{theorem1}, these observations immediately imply the following.

\begin{corollary}\label{corollary1}
Let $k$ be a positive integer.
It is NP-hard to decide whether $\gamma_R(G)-\gamma_{r2}(G)=k$ for a given connected $K_4$-free graph $G$. 
\end{corollary}
We proceed to our results concerning induced subgraphs.

\begin{theorem}\label{theorem2}
A graph $G$ satisfies $\gamma_{r2}(H)=\gamma_R(H)$ for every induced subgraph $H$ of $G$ if and only if $G$ is $\{P_5,C_5,C_4\}$-free. 
\end{theorem}
{\it Proof:} Since $\gamma_{r2}(H)<\gamma_R(H)$ for every graph $H$ in $\{P_5,C_5,C_4\}$, the necessity follows immediately.
In view of Theorem \ref{theorembefore}, 
in order to complete the proof, 
it suffices to show 
that every $\{P_5,C_5,C_4\}$-free graph $G$ satisfies $\gamma_R(G)\leq \gamma_{r2}(G)$.
Therefore, let $G$ be a $\{P_5,C_5,C_4\}$-free graph, 
and let $f$ be a minimum $2$-rainbow dominating function of $G$
such that $|f^{-1}(\{ 1,2\})|$ is as large as possible.
For $F\in 2^{\{ 1,2\}}$, let $V_{F}=f^{-1}(F)$.
Let
$$g:V(G)\to \{ 0,1,2\}:x\mapsto
\begin{cases}
0 & \mbox{, if $x\in V_{\emptyset}$,}\\
1 & \mbox{, if $x\in V_{\{ 1\}}\cup V_{\{ 2\}}$, and}\\
2 & \mbox{, if $x\in V_{\{ 1,2\}}$}.
\end{cases}
$$
Note that $w(g)=w(f)$.
If $g$ is a Roman dominating function of $G$, 
then $\gamma_R(G)\leq \gamma_{r2}(G)$.
Hence, 
we may assume that $g$ is not a Roman dominating function of $G$,
which implies the existence of a vertex $u$ in $V_{\emptyset}$ 
that has a neighbor $v_1$ in $V_{\{ 1\}}$ as well as a neighbor $v_2$ in $V_{\{ 2\}}$ 
but no neighbor in $V_{\{ 1,2\}}$.
We say that $v_1uv_2$ is a {\it special path}.

First, suppose that there is no special path $v_1uv_2$ such that $v_1$ and $v_2$ are adjacent,
that is, every special path is induced.
Let $v_1uv_2$ be a special path.
Let 
$$f_1:V(G)\to 2^{\{ 1,2\}}:x\mapsto
\begin{cases}
\emptyset &\mbox{, if $x\in \{ v_1,v_2\}$,}\\
\{ 1,2\} &\mbox{, if $x=u$, and}\\
f(x) & \mbox{, otherwise}
\end{cases}
$$
Since $w(f_1)=w(f)$, the choice of $f$ implies that $f_1$ is not a $2$-rainbow dominating function of $G$.
By symmetry, we may therefore assume that there is a special path $v_1u'v_2'$ with $u'\not\in N_G[u]$.
By our assumption, $v_1$ is not adjacent to $v_2$ or to $v_2'$.
Since $G$ is $C_4$-free, $u'$ is not adjacent to $v_2$, which implies that $v_2'$ is distinct from $v_2$.
Since $G$ is $C_4$-free, $u$ is not adjacent to $v'_2$.
Now, $G[\{ v_2',u',v_1,u,v_2\}]$ is $C_5$ or $P_5$ depending on whether $v_2$ and $v_2'$ are adjacent or not,
which is a contradiction.
Hence, there is a special path that is not induced.
If $v_1uv_2$ is a special path, and $v_1$ is adjacent to $v_2$, then we say that $v_1uv_2$ is a {\it special triangle}.

Let $U$ be a set of vertices of maximum order such that every vertex in $U$ belongs to some special triangle $T$ with $V(T)\subseteq U$.
Since there is at least one special triangle, the set $U$ is not empty.
For $F\in \{ \emptyset,\{ 1\},\{ 2\}\}$, let $U_F=U\cap V_F$.
By symmetry, we may assume that $|U_{\{ 1\}}|\geq |U_{\{ 2\}}|$.
Let 
$$f_2:V(G)\to 2^{\{ 1,2\}}:x\mapsto
\begin{cases}
\emptyset & \mbox{, if $x\in U_{\{ 1\}}$,}\\
\{ 1,2\} & \mbox{, if $x\in U_{\{ 2\}}$, and}\\
f(x) & \mbox{, otherwise.}
\end{cases}
$$
Since $w(f_2)\leq w(f)$, the choice of $f$ implies that $f_2$ is not a $2$-rainbow dominating function of $G$.
Since every vertex in $U_{\emptyset}\cup U_{\{ 1\}}$ has a neighbor in $U_{\{ 2\}}$, 
together with the definition of $U$, 
this implies the existence of 
a special triangle $v_1uv_2$ as well as a special path $v_1u'v_2'$ such that 
\begin{itemize}
\item $u\in U_{\emptyset}$, $v_1\in U_{\{ 1\}}$, and $v_1\in U_{\{ 2\}}$,
\item $u',v_2'\not\in U$, and 
\item $u'$ is not adjacent to $v_2$.
\end{itemize}
If $v_1$ and $v_2'$ are adjacent, then $v_1u'v_2'$ is a special triangle, 
and adding $u'$ and $v_2'$ to $U$ yields a contradiction to the choice of $U$.
Hence, $v_1$ is not adjacent to $v_2'$.
Since $G$ is $C_4$-free, $v_2$ is not adjacent to $v_2'$.
Let 
$$f_3:V(G)\to 2^{\{ 1,2\}}:x\mapsto
\begin{cases}
\emptyset & \mbox{, if $x=v_2$,}\\
\{ 1,2\} & \mbox{, if $x=v_1$, and}\\
f(x) & \mbox{, otherwise.}
\end{cases}
$$
Since $w(f_3)=w(f)$, the choice of $f$ implies that $f_3$ is not a $2$-rainbow dominating function of $G$.
This implies the existence of a vertex $u''\in V_{\emptyset}$
that is adjacent to $v_2$ but not to $v_1$.
Since $G$ is $C_4$-free, $u'$ is not adjacent to $u''$.
Now, $G[\{ v_2',u',v_1,v_2,u''\}]$ is $C_5$ or $P_5$ depending on whether $v'_2$ and $u''$ are adjacent or not,
which is a contradiction and completes the proof. $\Box$

\medskip

\noindent For a positive integer $k$, let ${\cal G}_k\left(\gamma_{R},\frac{3}{2}\gamma_{r2}\right)$
be the set of all graphs $G$ such that $\gamma_{R}(H)=\frac{3}{2}\gamma_{r2}(H)$ 
for every induced subgraph $H$ of $G$ with $\gamma_{r2}(H)\geq k$, that is,
$${\cal G}_k\left(\gamma_{R},\frac{3}{2}\gamma_{r2}\right)
=\left\{G:\forall H\subseteq_{\rm ind}G:\gamma_{r2}(H)\geq k\Rightarrow \gamma_{R}(H)=\frac{3}{2}\gamma_{r2}(H)\right\}.$$
Since $\gamma_{r2}(K_1)=1=\gamma_R(K_1)$,
the set ${\cal G}_1\left(\gamma_{R},\frac{3}{2}\gamma_{r2}\right)$ contains no graph of positive order.
Since $\gamma_{r2}(\bar{K}_2)=2=\gamma_R(\bar{K}_2)$,
the set ${\cal G}_2\left(\gamma_{R},\frac{3}{2}\gamma_{r2}\right)$ consists exactly of all complete graphs.  

\begin{theorem}\label{theorem3}
A graph $G$ belongs to ${\cal G}_3\left(\gamma_{R},\frac{3}{2}\gamma_{r2}\right)$ if and only if $G$ is $\{\bar{K}_3, K_2\cup K_1\}$-free.
\end{theorem}
{\it Proof:}
Since $\gamma_{r2}(\bar{K}_3)=\gamma_R(\bar{K}_3)=\gamma_{r2}(K_2\cup K_1)=\gamma_R(K_2\cup K_1)=3$, 
the graphs in ${\cal G}_3\left(\gamma_{R},\frac{3}{2}\gamma_{r2}\right)$ are $\{\bar{K}_3, K_2\cup K_1\}$-free.
In view of Theorem \ref{theorembefore}, 
in order to complete the proof, it suffices to show 
that every $\{\bar{K}_3, K_2\cup K_1\}$-free graph $G$ with $\gamma_{r2}(G)\geq 3$ 
satisfies $\gamma_{R}(G)=\frac{3}{2}\gamma_{r2}(G)$.
Therefore, let $G$ be a $\{\bar{K}_3, K_2\cup K_1\}$-free graph with $\gamma_{r2}(G)\geq 3$.
Since $\gamma_{r2}(G)\geq 3$, the graph $G$ is not complete.
Let $u$ and $v$ be two distinct vertices of $G$ that are not adjacent.
Since $\{\bar{K}_3, K_2\cup K_1\}$-free, we obtain $N_G(u)=N_G(v)=V(G)\setminus \{ u,v\}$.
This implies that
$$x\mapsto
\begin{cases}
\{ 1\} &\mbox{, $x=u$,}\\
\{ 2\} & \mbox{, $x=v$, and}\\
\emptyset & \mbox{, otherwise}
\end{cases}
$$
defines a $2$-rainbow dominating function of $G$ of weight $2$, which is a contradiction. $\Box$

\medskip

\noindent Theorem \ref{theorem1} implies that the graphs 
$G$ with $\gamma_{r2}(G)=\gamma_R(G)$ 
do not have a simple structure. 
In contrast to that, the graphs $G$ with 
$\gamma_R(G)= \frac{3}{2}\gamma_{r2}(G)$
seem far more restricted.
In fact, it is conceivable that these graphs 
can be recognized in polynomial time.
In our last result, we collect several of their properties.

\begin{theorem}\label{theorem4}
If $G$ is a graph with 
$\gamma_R(G)= \frac{3}{2}\gamma_{r2}(G)$, 
then every minimum $2$-rainbow dominating function $f$ of $G$
has the following properties,
where $V_F=f^{-1}(F)$ and $n_F=|V_F|$ for $F\in 2^{\{ 1,2\}}$.
\begin{enumerate}[(i)]
\item $n_{\{ 1\}}=n_{\{ 2\}}$ and $n_{\{ 1,2\}}=0$.
\item There are no edges between $V_{\{ 1\}}$ and $V_{\{ 2\}}$.
\item For $i\in [2]$, the maximum degree of $G[V_{\{ i\}}]$ is at most $1$.
\item For $i\in [2]$, every vertex in $V_{\emptyset}$ has at least $1$ and at most $2$ neighbors in $V_{\{ i\}}$.
\item For $i\in [2]$, every vertex $u$ in $V_{\{ i\}}$ has at least $2$ neighbors $v$ in $V_{\emptyset}$ 
with $N_G(v)\cap V_{\{ i\}}=\{ u\}$.
\end{enumerate}
\end{theorem}
{\it Proof:} Let $G$ be a graph with $\gamma_R(G)= \frac{3}{2}\gamma_{r2}(G)$, and let $f$ be a minimum $2$-rainbow dominating function of $G$.

\medskip

\noindent (i) Since the inequality (\ref{e1}) 
in the proof of Theorem \ref{theorembefore}
is satisfied with equality, 
we obtain $n_{\{ 1\}}=n_{\{ 2\}}$.
Similarly, since (\ref{e2}) is satisfied with equality, 
we obtain $n_{\{ 1,2\}}=0$.
Since $f$ is a $2$-rainbow dominating function of $G$ and $n_{\{ 1,2\}}=0$,
every vertex in $V_{\emptyset}$ has a neighbor in $V_{\{ 1\}}$ as well as a neighbor in $V_{\{ 2\}}$.

\medskip

\noindent (ii) If $v_1$ in $V_{\{ 1\}}$ is adjacent to $v_2$ in $V_{\{ 2\}}$, then 
$$
x\mapsto
\begin{cases}
0 & \mbox{, $x\in \{ v_2\}\cup V_{\emptyset}$,}\\
1 & \mbox{, $x\in V_{\{ 2\}}\setminus \{ v_2\}$, and}\\
2 & \mbox{, $x\in V_{\{ 1\}}$}
\end{cases}
$$
defines a Roman dominating function of $G$ of weight $\frac{3}{2}\gamma_{r2}(G)-1$, which is a contradiction.

\medskip

\noindent (iii) If $u$ in $V_{\{ i\}}$ is adjacent to two distinct vertices $v$ and $w$ in $V_{\{ i\}}$, then 
$$
x\mapsto
\begin{cases}
0 & \mbox{, $x\in \{ v,w\}\cup V_{\emptyset}$,}\\
1 & \mbox{, $x\in V_{\{ i\}}\setminus \{ u,v,w\}$, and}\\
2 & \mbox{, $x\in \{ u\}\cup V_{\{ 3-i\}}$}
\end{cases}
$$
defines a Roman dominating function of $G$ of weight $\frac{3}{2}\gamma_{r2}(G)-1$, which is a contradiction.

\medskip

\noindent (iv) As observed above, every vertex in $V_{\emptyset}$ has a neighbor in $V_{\{ i\}}$.
If $u$ in $V_{\emptyset}$ is adjacent to three distinct vertices $v_1$, $v_2$, and $v_3$ in $V_{\{ i\}}$, then 
$$
x\mapsto
\begin{cases}
0 & \mbox{, $x\in \{ v_1,v_2,v_3\}\cup (V_{\emptyset}\setminus \{ u\})$,}\\
1 & \mbox{, $x\in V_{\{ i\}}\setminus \{ v_1,v_2,v_3\}$, and}\\
2 & \mbox{, $x\in \{ u\}\cup V_{\{ 3-i\}}$}
\end{cases}
$$
defines a Roman dominating function of $G$ of weight $\frac{3}{2}\gamma_{r2}(G)-1$, which is a contradiction.

\medskip

\noindent (v) Let $i\in [2]$ and let $u\in V_{\{ i\}}$.
Let $P(u)=\{ v\in V_{\emptyset}: N_G(v)\cap V_{\{ i\}}=\{ u\}\}$.
If $P(u)=\emptyset$, then
$$
x\mapsto
\begin{cases}
0 & \mbox{, $x\in V_{\emptyset}$,}\\
1 & \mbox{, $x\in V_{\{ 3-i\}}\cup\{ u\}$, and}\\
2 & \mbox{, $x\in V_{\{ i\}}\setminus \{ u\}$}
\end{cases}
$$
defines a Roman dominating function of $G$ of weight $\frac{3}{2}\gamma_{r2}(G)-1$, which is a contradiction.
Hence, $P(u)$ is non-empty.
If $P(u)=\{ v\}$, then let $w$ be a neighbor of $v$ in $V_{\{ 3-i\}}$.
Now,
$$
x\mapsto
\begin{cases}
0 & \mbox{, $x\in \{ u,w\}\cup (V_{\emptyset}\setminus \{ v\})$,}\\
1 & \mbox{, $x\in V_{\{ 3-i\}}\setminus \{ w\}$, and}\\
2 & \mbox{, $x\in \{ v\}\cup (V_{\{ i\}}\setminus \{ u\})$}
\end{cases}
$$
defines a Roman dominating function of $G$ of weight $\frac{3}{2}\gamma_{r2}(G)-1$, which is a contradiction.
$\Box$

\medskip

\noindent {\bf Acknowledgment}  
J.D. Alvarado and S. Dantas were partially supported by FAPERJ, CNPq, and CAPES.
D. Rautenbach was partially supported by CAPES.

\end{document}